\documentclass{article}
\usepackage{amsmath}
\usepackage{amsfonts}
\usepackage{amssymb}	
\usepackage{color}
\newtheorem{theo}{Theorem}[section]

\newtheorem{rem}[theo]{Remark}
\newtheorem{prop}[theo]{Proposition}
\newtheorem{defi}{Definition}[section]

\newcommand{\mysection}[1]{\section{#1} \setcounter{equation}{0}}
\newcommand{\proof}{{\sc Proof.} \quad}

\newcommand{\R}{\mathbb{R}}

\newcommand{\be}{\begin{equation} \label}
\newcommand{\ee}{\end{equation}}
\newcommand{\bes}{\begin{equation} \begin{array}{c} \label}
\newcommand{\ees}{\end{array} \end{equation}}
\newcommand{\bea}{\begin{eqnarray}\label}
\newcommand{\eea}{\end{eqnarray}}
\newcommand{\beas}{\begin{eqnarray} \begin{array}{rcl} \label}
\newcommand{\eeas}{\end{array} \end{eqnarray}}
\newcommand{\bas}{\begin{eqnarray*}}\newcommand{\eas}{\end{eqnarray*}}
\newcommand{\bass}{\begin{eqnarray*} \begin{array}{rcl}}
\newcommand{\eass}{\end{array} \end{eqnarray*}}
\newcommand{\basss}{\begin{eqnarray*} \begin{array}{c}}
\newcommand{\easss}{\end{array} \end{eqnarray*}}
\newcommand{\qed}{{}\hfill $\square$ \\}
\newcommand{\bit}{\begin{itemize}}
\newcommand{\eit}{\end{itemize}}
\newcommand{\nn}{\nonumber}

\newcommand{\diam}{\mathrm{diam}\,}
\begin{document}
\title{A theorem on measures in dimension $2$ and applications to vortex sheets.}
\author{
Tomasz Cie\'slak \\
{\small Institute of Mathematics, Polish Academy of Sciences, \'Sniadeckich 8, Warsaw, Poland}\\
{\small e-mail: T.Cieslak@impan.pl}\\
Marta Szuma\'nska\\
{\small Institute of Mathematics, University of Warsaw, Banacha 2, Warsaw, Poland}\\
{\small e-mail: M.Szumanska@mimuw.edu.pl}}

\maketitle

\begin{abstract} We find conditions under which measures belong to $H^{-1}(\R^2)$. Next we show that measures generated by Prandtl, Kaden as well as Pullin spirals, objects considered by physicists as incompressible flows generating vorticity, satisfy assumptions of our theorem, thus they are (locally) elements of $H^{-1}(\R^2)$. Moreover, as a by-product, we prove an embedding of the space of Morrey type measures in $H^{-1}$.   

\end{abstract}

\mysection{Introduction}\label{section1}

Due to d'Alembert's paradox we know that in the class of regular steady irrotational Euler flows the lift exerted by the inviscid incompressible flow on the three-dimensional body is zero, see \cite[Appendix 1.4]{mar_pul}. On the other hand, among regular solutions to the Euler equation, initially potential flow stays potential during the evolution according to the Helmholtz theorem. Hence, physicists (not only) were looking for the irregular flows generating the vorticity.

It was Felix Klein, in his famous Kaffel\"offel experiment, postulating that the flow in which all the vorticity is supported on a spiral may generate the vorticity from an initially potential flow, \cite{klein}, \cite[Chapter 6.1]{saffman}.  Such spirals are examples of {\it vortex sheets} i.e. flows whose vorticity is supported on a curve, being zero off the curve. Later experiments (including wind tunnel experiments) were also suggesting the vortex sheet flows to be responsible for the creation of a vorticity from initially potential flow. Among theoretical and experimental examples of flows leading to the creation of vorticity were  creation of a circulation about a body, leading edge suction problem, development of the sheet into a pair of vortices, accelerated flow past a wing or flow separation,  see \cite[Chapter 6]{saffman}. 

The main theoretical efforts in this area were made in the school of Prandtl. The examples of
vortex sheets 
best known in the literature come from Prandtl, see \cite{prandtl}, Kaden, see \cite{kaden} or Anton, see \cite{anton}. All those examples are self-similar spirals of vorticity (vorticity off the spiral is zero). Broader class of such self-similar spirals come from the numerical experiments carried by Pullin, \cite{pullin}.
Such spirals appear as a result of experiments and theoretical considerations, see Section {\ref{section3}}, but it is still not known if, and in what sense, such objects can be interpreted as solutions to the 2d incompressible inviscid Euler equation. One has to keep in mind that even if some examples of flows generating the vorticity do not solve the 2d Euler equation they can still bring some important information, being for instance solutions to Euler-like problems. Well known example of such a phenomenon is the Prandtl-Munk steady vortex sheet. It has been introduced in the literature as a flow generating the vorticity, some physicists and specialists in computational fluid mechanics were considering it as a potential solution of 2d Euler equation. However, as it was shown in \cite{lop_nus_sou} it is actually a weak solution to the nonhomogeneous Euler equation with forcing term representing a tension force applied to the tips of a wing. 

The present paper is the first step of our investigations whose goal is understanding if known examples of flows generating the vorticity are weak solutions of the Euler system, or, if the answer is negative, if they are solutions to some Euler-like problems. Another issue of our long-term plan is understanding whether the known examples of flows generating the vorticity can be helpful in finding counterexamples to the uniqueness of Delort's solutions, see \cite{delort}. In the present paper we are mostly concerned with the so-called Prandtl's, Kaden's and Pullin's spirals.  
We give an answer to a question whether such objects are locally in $H^{-1}(\R^2)$, i.e. whether their restriction to any compact set belongs to $H^{-1}(\R^2)$. 
There are several reasons for studying the above question. 
  The framework of weak solutions to the 2d Euler equations capturing the vortex sheets introduced in the pioneering work \cite{DiP_Majda} consists of a few requirements concerning the velocity field $v$ in order to interpret it as a vortex sheet solution to the 2d Euler equations.  One of them is the following 
\begin{equation}\label{energy}
v\in L^\infty(0,T;L^2_{loc}(\R^2)).
\end{equation} 
It is not only the technical assumption so that weak velocity formulation makes sense. The physical interpretation of \eqref{energy} is clear, kinetic energy of the fluid is finite at least locally. Our question if the measure supported on the spiral of vorticity belongs locally to $H^{-1}$ is clearly related to \eqref{energy}. 

Delort proved a theorem stating the global existence of vortex sheet type weak solutions to the 2d Euler equation, (\cite{delort}). The assumptions he imposes are that the initial vorticity is a compactly supported Radon measure with a fixed sign (positive or negative) being an element of $H^{-1}$. Hence if one wants to see whether flows generating vorticity can be understood as Delort's solutions, one has to check if they initially belong to $H^{-1}$. Delort's theorem concerns only compactly supported data, however one may hope to extend the theorem to the non-compact case. 

Examples of compactly supported spirals of vorticity are constructed by considering complex potentials, as it is in the case of the Anton spiral (see Section \ref{section3}), or by numerical experiments \cite{lop_low_nus_zheng}.
Both methods give a description  of spirals which seems to be not sufficient for further mathematical analysis.  
However, asymptotical behaviour of Anton's spiral around its center is similar to the behaviour of Kaden's spiral, thus
examining local properties of Kaden's or Prandtl's spirals may give an idea for features of known and future examples of compactly supported spirals of vorticity. The measure theoretic fact that can be useful in checking whether a  compactly supported spiral of vorticity is an element of $H^{-1}(\R^2)$ is the content of our main theorem (Theorem \ref{Tw.1}).

The numerical example from \cite{lop_low_nus_zheng} was a step towards constructing an evidence for nonuniqueness of Delort's solutions. Finding whether the solution is not unique is a long-standing open problem and Kaden type spirals (if one could obtain compactly supported ones) could have been a proper counterexample. Notice that one needs to check whether such a solution is initially in $H^{-1}$ to close the argument. Theorem \ref{Tw.1} may be a useful tool for this purpose. 

Let us also mention that the requirement that solution stays in $H^{-1}$ for times $t>0$ is very important for several reasons. First it guarantees \eqref{energy}, a necessary condition of being a weak solution of 2d Euler equation. Second the fact that solution stays in $H^{-1}$ excludes the finite-time collapse into a point vortex. Moreover, it guarantees the equivalence of weak velocity and weak vorticity formulations of the 2d Euler equation, see \cite{schochet}. Finally, if one does not require any regularity of vorticity, and is fully satisfied with solutions satisfying weak velocity formulations, then vortex sheet evolution is not unique according to the counterexamples coming from the theory of wild solutions, see \cite{szekelyhidi}. 

Another reason for considering whether spirals of vorticity are elements of $H^{-1}$ is related to 
the Birkhoff-Rott equation. We will introduce the Birkhoff-Rott equations in details in Section \ref{section1.5}, now let us only mention that the Birkhoff-Rott equation is a description of the evolution of the interface between the irrotational parts of a flow, which is a self-similar spiral in our case. It is important to know whether the two descritpions, by the weak solutions of Euler's equation or solution to the Birkhoff-Rott equation give the same result. Since the self-similar spirals of vorticity are not smooth curves, one cannot apply standard theory (presented in \cite{mar_pul}) and has to use the recent results shown in \cite{lop_nus_sch}. It is proved there that if a weak solution to either Birkhoff-Rott or Euler equations stays in $H^{-1}$ when time evolves, then both descriptions of the evolution are equivalent.    
            
Our proof of the fact that the Prandtl and Kaden spirals are locally in $H^{-1}(\R^2)$ is based on a theorem which is a more general tool for recognizing whether compactly supported Radon measures are elements of $H^{-1}(\R^2)$.
\begin{theo}\label{Tw.1}
Let $\mu$ be a positive Radon measure with a support inside a ball $B(0,R_0)\subset \R^2$. Assume that  there exists a positive constant $c_1$ such that for any $r<R_0$
\begin{equation}\label{warunek}
\mu(B(0,r))=c_1 r^{\alpha},\;\;\mbox{where}\;\;\alpha>0, 
\end{equation}
then $\mu\in H^{-1}(\R^2)$.
\end{theo}

Notice that due to the linear structure of $H^{-1}$ and the Hahn decomposition theorem one can use Theorem \ref{Tw.1}  to deal with signed Radon measures. In particular the following remark holds.  
\begin{rem}\label{uwaga1}
Let $\mu_+$ and $\mu_-$ be respectively positive and negative part of a measure with a support included in the ball $B(0,R_0)\subset \R^2$. If both $\mu_+$ and $\mu_-$ satisfy \eqref{warunek} with exponents $\alpha_+$ and $\alpha_-$ for $r<R_0$, possibly $\alpha_+\neq\alpha_-$, then $\mu\in H^{-1}(\R^2)$.
\end{rem} 
Moreover, we state the following remark, which asserts the applicability of our results for spirals with an origin in $0\neq x_0\in \R^2$. 

\begin{rem}\label{uwaga2.0}
In \eqref{warunek} the origin of the mass can be shifted to any fixed point $x_0$, it does not influence the finiteness of the norm in $H^{-1}(\R^2)$.
\end{rem}

Let us notice that condition \eqref{warunek} is natural in applications to vortex sheets as it is satisfied by self-similar spirals of vorticity undergoing the laws of similitude, see Section \ref{section3}. 

Last but not least, the methods we use in the proof of Theorem \ref{Tw.1} yield a direct proof of the continuous embedding of spaces of Morrey type measures in $H^{-1}(\R^n)$, $n\geq 1$, what is one of the claims of \cite[Theorem 4.3]{lop_nus_tad}. 

\mysection{Preliminaries on vortex sheets.}\label{section1.5}

As it was mentioned in the introduction a vortex sheet corresponds to the situation when the whole vorticity of a velocity field is concentrated on a curve. To avoid ambiguity we specify  that by a plane curve we understand a subset of $\R^2$ (or a complex plane) which is a continuous and injective image of a circle or of an (possibly open) interval $I\subset \R$.  For convenience we will not distinguish a curve and its parametrization. 

A vortex sheet is not only a curve which describes its shape, but also a function describing a strength of vorticity given at each point that will be denoted by $\gamma$. The function can be considered as a density of a measure supported on a curve, thus mathematically vortex sheet represents and is represented by a certain type of singular measures on the plane.  Before we formulate propositions that will clarify the physical meaning of $\gamma$ let us specify our general setting. 

We assume that we are given a plane curve $\eta$ and a velocity field $v$ which is $C^1$ on  $\R^2 \setminus \eta$.
We assume also that the field is incompressible (i.e. divergence free) in the sense of distributions and irrotational off the curve $\eta$. 
Notice that, by the definition, curve $\eta$ divides any, sufficiently small, ball centered on $\eta$ into two disjoint parts; as we may assume that the curve is oriented we can distinguish the left and the right part of such ball. We will denote by $v^+$ the velocity field on the left part, while by $v^-$ on the right part. If curve $\eta$ cuts a plane in two disjoint parts, then the vector field $v$ on $\R^2\setminus \eta$ can be divided into two vector fields $v^+$ (on the left side of the curve) and $v^-$. 

Using the above assumptions and notation we can state the following two propositions.

\begin{prop}\label{prop3.0.1}
Let $\eta$ be a $C^1$ plane curve. Assume that a velocity field $v$ on $\R^2$ is $C^1$ besides a curve $\eta$, it is also divergence-free in the sense of distributions. Then, there holds $(v^+-v^-)\cdot\vec{n}=0$ on the curve $\eta$ , where $\vec{n}$ is a normal vector to the curve $\eta$.
\end{prop}
\proof
The claim follows from the distributional incompressibility condition 
and the Green formula.
  
\begin{prop}\label{prop3.0.2}
Let $\eta$ be a $C^1$ plane curve. Assume that a velocity field $v$ on $\R^2$ to be $C^1$ besides a curve $\eta$ and irrotational off the curve. Then rot v equals, in the sense of distributions, to a measure supported on $\eta$ with a density $\gamma$, $\gamma:=(v^+-v^-)\cdot\tau$, where $\tau$ denotes a tangent vector to the curve $\eta$.
\end{prop}  
\proof
For each poit of $\eta$ we fix a ball centered at this point such that $\eta$ cuts the ball into two disjoint parts.
By $\Omega$ we denote the sum of those balls. As curve $\gamma$ cuts $\Omega$ into two disjoint part we can distinguish the left part $\Omega_+$ and the right $\Omega_-$ equipped with the velocity field $v^+$ and $v^-$ respectively.

Take a smooth function $\phi$ compactly supported in $\Omega\in \R^2$. 
We compute a vorticity of $v$ in the sense of distributions. By the Green formula
\[ 
\int_\Omega rot v\phi dx=\int_{\Omega_+}rot(v^+\phi)dx+ \int_{\Omega_-}rot(v^-\phi)dx=
\int_{\eta\cap supp_{\phi}}\phi(v^+-v^-)\cdot\tau ds.
\] 
\qed

In Section \ref{section3} we consider curves which are not  $C^1$, but that are at least $C^1$ regular besides an origin, being continuous up to an origin. 
Notice that in general we do not assume the curve to be locally rectifiable at the origin, but we guarantee that a measure defined by the strength of vorticity is  finite on every ball. Then, by \cite[Theorem 2.18]{rudin}, we know that such an object is a Radon measure.

For any $z\in \eta$ by $\eta_z$ we denote a part of $\eta$ connecting $0$ and $z$. Cumulated vorticity $\Gamma$ at $z$ is given by the formula
\begin{equation}\label{Gam}
\Gamma(z):=\int_{\eta_z}\gamma(z).
\end{equation}
One immediately sees that for a measure 
 supported on a curve $\eta$ with a density $\gamma$, $\Gamma(z)$ is equal to the
measure of arc $\eta_z$.  As it should not lead to cofusion, we will use the same letter to denote the measure on $\R^2$ generated by $\Gamma$ i.e. for  a subarc $\eta_s \subset \eta$ connecting two points $z_1$ and $z_2$, we write $\Gamma(\eta_s):= \Gamma(z_2) - \Gamma(z_1)$ and for any Borel set $A\subset \R^2$ we put $\Gamma(A) := \Gamma(A\cap \eta)$.

Notice also that if the strength of vorticity has a constant sign or the curve can be divided into two parts on which of each vorticity has a constant sign, then $\Gamma(z)$ is an injective function defined on $\eta$ and the value of $\Gamma$ uniquely defines a point on the curve, thus one can use $\Gamma$ to parametrize $\eta$. 

The Birkhoff-Rott equation describes the evolution of a vortex sheet. 
The position of the curve at the time $t$ and cumulated vorticity $\Gamma$ is desribed by the complex variable $z(\Gamma, t)$ with the use of the following equation.
\begin{equation}\label{B-R}
\frac{d}{dt}\bar{z}(\Gamma, t)=\frac{1}{2\pi i}p.v.\int\frac{d\Gamma'}{z(\Gamma,t)-z(\Gamma',t)}.
\end{equation}
The equivalence between the description of the vortex sheet through the Birkhoff-Rott equation and the 2d Euler equation is studied in \cite{mar_pul} for the smooth interface and in \cite{lop_nus_sch} for less regular curves. 

\mysection{Measure theoretic results. The proof of the main theorem.}\label{section2}

 Before we proceed with the proof of Theorem \ref{Tw.1} let us recall several necessary definitions. The first one, see \cite{rauch}, in the special case $s=-1$ specifies our target space $H^{-1}$.

\begin{defi}\label{defi0}
A tempered distribution $f$ belongs to $H^s(\R^n)$, $s\in\R, n\geq 1$ if
\[
\left\|f\right\|_{H^{-1}(\R^n)}^2:=\int_{\R^n}\left(1+|x|^2\right)^s|\hat{f}(x)|^2dx<\infty,
\]
where $\hat{f}$ denotes the Fourier transform of $f$.
\end{defi}
Crucial tools used in the proof of Theorem \ref{Tw.1} are t-energies of a positive Radon measure $\mu$ (see \cite{mattila}), which are helpful when one wants to control the concentration of $\mu$ on small balls.
\begin{defi}\label{defi1}
For a positive measure $\mu$ on $\R^n$, $n\geq 1$, we define for any $t>0$ t-energy related to the measure $\mu$ by
\begin{equation}\label{t-energy}
I_t(\mu)=\int_{\R^n}\int_{\R^n}|x-y|^{-t}d\mu (x)d\mu(y).
\end{equation} 
\end{defi}
For $0<t<n$ t-energy $I_t(\mu)$ can be expressed in terms of a Fourier transform of a compactly supported positive Radon measure, see \cite[Lemma 12.12]{mattila}, namely there exists a positive constant $c(t,n)$ such that
\begin{equation}\label{Is}
I_t(\mu)=(2\pi)^{-n}c(t,n)\int_{\R^n}|x|^{t-n}|\hat{\mu}(x)|^2dx.
\end{equation}
Next notice that in view of (\ref{Is}) and the fact that for any $s>0$ and all $x\in \R^n$
\[
\left(1+|x|^2\right)^{-1}\leq |x|^{s-2},
\]
the following proposition holds.
\begin{prop}\label{prop2.1}
Let $\mu$ be a positive Radon measure with compact support in $\R^2$. If there exists $0<s<2$ such that $I_s(\mu)<\infty$ then $\mu\in H^{-1}(\R^2)$.
\end{prop}
One of the components of the proof of Theorem \ref{Tw.1} are concentration estimates being a consequence of  \eqref{warunek}. We state them in a proposition below. 
\begin{prop}\label{prop2.2}
Assume that a positive Radon measure with a support inside $B(0,R_0)\subset\R^2$ satisfies \eqref{warunek}. Then for any $\{x: 2r<|x|\leq R_0\}$ the following conditions are satisfied 
\begin{equation}\label{xkulki}
\mu(B(x,r))\leq C(c_1,\alpha,R_0)\; r\;\;\mbox{if}\;\;\alpha\geq 1  
\end{equation}
and
\begin{equation}\label{xkulki1}
\mu(B(x,r))\leq C(c_1,\alpha)\; r|x|^{\alpha-1}\;\;\mbox{for}\;\;0<\alpha<1.
\end{equation} 
\end{prop}  
Before starting the proof of Proposition \ref{prop2.2} let us
notice that \eqref{xkulki1} does not contradict the local finiteness of a measure $\mu$ at the origin. Namely, when $x\rightarrow 0$ then $\mu\left(B(x,r)\right)\leq C|x|^{\alpha-1}r\leq C r^{\alpha}$. 
\proof
For $\{(x,r):2r< |x|\leq R_0\}$ it is clear that 
\[
B(x,r)\subset B(0,|x|+r)\setminus B(0,|x|-r).
\]
Hence 
\[
\mu(B(x,r))\leq C(c_1,\alpha)|x|^{\alpha}\left[\left(1+\frac{r}{|x|}\right)^\alpha-\left(1-\frac{r}{|x|}\right)^\alpha\right].
\]
The Taylor expansion of the right-hand side of the above formula up to the second order terms gives
\[
\mu(B(x,r))\leq C(c_1,\alpha)|x|^{\alpha}\left[1+\alpha\frac{r}{|x|}-1+\alpha\frac{r}{|x|}+\frac{(\alpha|\alpha-1|)}{2}\left((1+d_1)^{\alpha-2}+(1-d_2)^{\alpha-2}\right)\frac{r^2}{|x|^2}\right],
\]
where $d_1,d_2\in(0,\frac{r}{|x|})\subset(0,1/2)$. Hence
\begin{eqnarray}\label{2.0}
\mu(B(x,r))&\leq& C(c_1,\alpha)|x|^{\alpha}\left(2\alpha\frac{r}{|x|}+c(\alpha)\frac{r^2}{|x|^2}\right)\nn \\
&\leq& C(c_1,\alpha)\left(2\alpha r|x|^{\alpha-1}+c(\alpha)r^2|x|^{\alpha-2}\right).
\end{eqnarray} 
Thus, since $r<\frac{|x|}{2}$, for $\alpha\geq 1$ in view of \eqref{2.0} we have,
\[
\mu(B(x,r))\leq C(c_1,\alpha)r|x|^{\alpha-1}\leq C(c_1,\alpha)R_0^{\alpha-1}r
\]
and \eqref{xkulki} follows. For $\alpha<1$ on the other hand, again using $\frac{r}{|x|}<\frac{1}{2}$ in \eqref{2.0}, we arrive at
\[
\mu(B(x,r))\leq C(c_1,\alpha)r|x|^{\alpha-1}+c(\alpha)r^2|x|^{\alpha-2}\leq C(c_1,\alpha)r|x|^{\alpha-1}
\]
and \eqref{xkulki1} follows.
\qed

Now we are ready to proceed with the proof of Theorem \ref{Tw.1}.

{\bf Proof of Theorem \ref{Tw.1}.}

The proof consists of two parts. In the first one we deal with an easier case when $\alpha\geq 1$, in the second we prove the claim in the remaining range of parameters $\alpha\in (0,1)$, when the measure of a ball depends on the position of its origin. 
\begin{enumerate}
  \item[$\alpha\geq 1$ ] Using the Fubini theorem and the change of variables $r=u^{-\frac{1}{s}}$ we obtain
\begin{eqnarray}\label{2.1}
I_s(\mu)&=&\int_{\R^2}\int_{\R^2}|x-y|^{-s}d\mu(y)d\mu(x)
=\int_{\R^2}\int_0^\infty \mu(\{y:|x-y|^{-s}\geq u\})dud\mu(x)\nn \\
&=&\int_{\R^2}\int_0^\infty \mu \left(B\left(x,u^{-\frac{1}{s}}\right)\right)du d\mu(x)=\int_{\R^2}s\int_0^\infty r^{-s-1}\mu(B(x,r))drd\mu(x)\nn \\
&\stackrel{(\ref{xkulki})}\leq &C(c_1,\alpha,R_0)s\int_{B(0,R_0)}\int_0^{\frac{|x|}{2}}r^{-s}drd\mu(x)\nn \\
&+ & \int_{B(0,R_0)}s\int_{\frac{|x|}{2}}^\infty r^{-s-1}\mu(B(0,R_0))drd\mu(x),
\end{eqnarray}
We see that in order to make sure that the last two terms are integrable, we have to restrict $0<s<1$. Then from \eqref{2.1} we infer
\[
I_s(\mu)\leq C(c_1,\alpha,R_0)s\int_{B(0,R_0)}\frac{|x|^{1-s}}{(1-s)2^{1-s}}d\mu(x)+2^s\mu(B(0,R_0))\int_{B(0,R_0)}|x|^{-s} d\mu(x).
\]
Since $\mu$ is positive and compactly supported, we observe that  
\begin{eqnarray}\label{fin_arg}
\int_{B(0,R_0)} |x|^{-s}d\mu(x)&< &\int_{\R^2}|x|^{-s}d\mu(x)\nn \\ 
&= &\int_0^\infty\mu(\{x:|x|^{-s}\geq u\})du
=s\int_0^\infty r^{-s-1}\mu \left(B(0,r)\right)dr\nn \\
&\leq &s\left(c_1\int_0^{R_0} r^{-s-1}r^\alpha dr+\int_{R_0}^\infty r^{-s-1}\mu(B(0,R_0))dr\right).
\end{eqnarray}
The first integral in \eqref{fin_arg} is finite if $\alpha-s>0$, the finiteness of the second one is provided by $s>0$. Hence for $s\in(0,\min \{\alpha,1\})=(0,1)$ we see that $I_s<\infty$ and so, by Proposition \ref{prop2.1}, $\mu\in H^{-1}(\R^2)$.

	\item[$\alpha \in (0,1)$]The scheme of the reasoning is the same as in the previous case. Again we start with computing $I_s(\mu)$ in order to infer the claim from Proposition \ref{prop2.1}. As in \eqref{2.1} we use Fubini's theorem and the same change of variables. 
\begin{eqnarray}\label{2.2}
I_s(\mu)&=&\int_{\R^2}\int_{\R^2}|x-y|^{-s}d\mu(y)d\mu(x)= \int_{\R^2}s\int_0^\infty r^{-s-1}\mu\left(B(x,r) \right)drd\mu(x)\nn \\
&\leq& C(c_1,\alpha)s \int_{\R^2}\int_0^{\frac{|x|}{2}}r^{-s}|x|^{\alpha-1}drd\mu(x)+ \int_{\R^2}s\int_{\frac{|x|}{2}}^\infty r^{-s-1}\mu(B(0,R_0))drd\mu(x)\nn \\
&=&C(c_1,\alpha)s \frac{1}{(1-s)2^{1-s}}\int_{\R^2} |x|^{\alpha-s}d\mu(x)+2^s\mu(B(0,R_0))\int_{\R^2} |x|^{-s}d\mu(x). 
\end{eqnarray}
First we observe that in the second line the term $\int_0^{\frac{|x|}{2}}r^{-s}|x|^{\alpha-1}dr$ is finite since  we assume $s<1$. Assuming additionally $0<s<\alpha<1$ we show the finiteness of both terms on the right-hand side of \eqref{2.2}. Indeed, in the first term 
\[
\int_{\R^2} |x|^{\alpha-s}d\mu(x)=\int_{B(0,R_0)} |x|^{\alpha-s}d\mu(x)\leq R_0^{\alpha-s}\mu(B(0,R_0)).
\] 
The second term has already been estimated in \eqref{fin_arg}. Therefore $I_s(\mu)<\infty$ for $0<s<\alpha$ and 
$\mu\in H^{-1}(\R^2)$.  
 
\end{enumerate} 
\qed

Let us also emphasize that using the t-energy methods we obtain a direct measure theoretic proof of continuous embeddings of spaces of Morrey type measures in $H^{-1}$. Before we state the embedding theorem and present the proof, we give the necessary definitions following \cite{lop_nus_tad}. 
\begin{defi}\label{defi2}
We say that $\mu\in M^{p}(\R^n)$, the subspace of Radon measures, if 
\begin{equation}\label{Morrey}
\left\|\mu\right\|_{M^p}:=\sup_{R>0}\left[R^{-n\left(1-\frac{1}{p}\right)}\sup_{x\in \R^n}|\mu|(B(x,R))\right]<\infty,
\end{equation}
where $|\mu|$ denotes the total variation of a measure.
\end{defi}
Let us recall that by the Hahn decomposition theorem every Radon measure $\mu$ can be decomposed into its positive and negative parts, i.e. there exist positive Radon measures $\mu^{+},\mu^{-}$ such that
\[
\mu=\mu_+-\mu_-.\;\;\mbox{Then}\;\;|\mu|=\mu_+ +\mu_-,
\] 
\begin{equation}\label{2.4}
\;\;\max\{\mu^{+}(A),\mu^{-}(A)\}\leq |\mu|(A),
\end{equation}
and
\begin{equation}\label{H-1}
\left\|\mu\right\|_{H^{-1}}\leq \left\|\mu_+\right\|_{H^{-1}}+\left\|\mu_-\right\|_{H^{-1}}.
\end{equation}
We begin with stating an $n$-dimensional analogue of Proposition \ref{prop2.1}.
\begin{prop}\label{prop2.5}
Let $\mu$ be a positive Radon measure on $\R^n$. Assume that $I_s(\mu)<\infty$ for $n>s>n-2$, then $\mu\in H^{-1}(\R^n)$. 
\end{prop}
\proof Notice that for $-2<s-n<0$ and all $x\in\R^n$
\[
(1+|\xi|^2)^{-1}<|\xi|^{s-n},
\]  
hence by \eqref{Is} we get
\[
\left\|\mu\right\|_{H^{-1}}^2\leq \int_{\R^n} \left|\hat{\mu}\right|^2|\xi|^{s-n}d\xi =(2\pi)^nc(t,n)^{-1}I_s(\mu).
\]
\qed
 
In the following text by $B{\cal M}$ we denote the space of bounded Radon measures and by $M^p(\Omega)$ we mean a subspace of $M^p(\R^n)$ consisting of measures supported in $\Omega$.
\begin{theo}\label{Tw.2}
Let $\Omega$ be a bounded subset of $\R^n$. Then $M^{p}(\Omega)\cap B{\cal M}$ is continuously embedded in $H^{-1}(\R^n)$ for $p>\frac{n}{2}$. 
\end{theo}
\proof
First notice that by \eqref{2.4}, \eqref{H-1} and Proposition \ref{prop2.5} 
\[
\left\|\mu\right\|_{H^{-1}}^2\leq I(\mu_+)+I(\mu_-)\leq 2 I_s(|\mu|).
\]
Similarly as in the proof of Theorem \ref{Tw.1}, Fubini's theorem and a change of variables lead to an inequality 
\begin{eqnarray}\label{mor}
\left\|\mu\right\|_{H^{-1}}^2&\leq& 2\int_{\R^n}\int_{\R^n}|x-y|^{-s}d|\mu|(y)d|\mu|(x)=2\int_{\R^n}s\int_0^\infty r^{-s-1}|\mu|(B(x,r))dr d|\mu|(x)\nn \\
&\leq& 2s\int_{\R^n}\left(\int_0^1 r^{-s-1}r^{n\left(1-\frac{1}{p}\right)}\left\|\mu\right\|_{M^p}dr+\int_1^\infty r^{-s-1} |\mu|(\Omega)dr\right)d|\mu|(x),
\end{eqnarray}
the bound on $|\mu|(B(x,r))$ in the second line is a consequence of \eqref{Morrey}. 

Moreover, we have
\begin{equation}\label{pp}
\left\|\mu\right\|_{M^p}\geq(\diam \Omega)^{-n\left(1-\frac{1}{p}\right)}\sup_{x\in \R^n}|\mu|(B(x,\diam \Omega))\geq (\diam \Omega)^{-n\left(1-\frac{1}{p}\right)}|\mu|(\Omega).
\end{equation}
Plugging the above estimate in \eqref{mor} we obtain
\begin{eqnarray}
\left\|\mu\right\|_{H^{-1}}^2&\leq & 2s\left\|\mu\right\|_{M^p}\int_{\R^n}\int_0^1r^{-s-1+n\left(1-\frac{1}{p}\right)}drd|\mu|\nn \\
&+ &2s (\diam \Omega)^{n\left(1-\frac{1}{p}\right)}\left\|\mu\right\|_{M^p}\int_{\R^n}\int_1^\infty r^{-s-1}drd|\mu|.\nn
\end{eqnarray}
Hence we see that for $0<s<n-\frac{n}{p}$, using once more \eqref{pp} we can continue and obtain
\[
\left\|\mu\right\|_{H^{-1}}^2\leq \left\|\mu\right\|_{M^p} C(\diam \Omega, s, n, p)\int_{\R^n}d|\mu|(x)\leq C(\diam \Omega, s, n, p)\left\|\mu\right\|_{M^p}^2.
\]
\qed
Let us comment that Theorem \ref{Tw.2} leads to an alternative prove of 
Theorem \ref{Tw.1} for the restricted range of parameters $\alpha\geq1$ (see the remark below). 
However, it seems to us that one cannot handle the case $\alpha<1$ this way, since the dependence on the position of the origin of the ball in the estimate \eqref{xkulki1} disturbs the possibility of applying Theorem \ref{Tw.2} 
 As we shall see in the next section, $\alpha<1$ corresponds to Kaden spirals which seem to be of importance from the point of view of physics as well as mathematics.

\begin{rem}\label{Kaden}
If a measure $\mu$ satisfies assumptions of Theorem \ref{Tw.1} for $\alpha\geq 1$, then by \eqref{xkulki} we find out that $\mu\in M^p(B(0,R))$ for any $p\in(1,2)$, therefore Theorem \ref{Tw.2} implies that $\mu\in H^{-1}(\R^2)$.
\end{rem}

\mysection{Applications to spirals of vorticity}\label{section3}

In the current section we present the spirals of vorticity introduced by physicists. Our aim is to show that the Radon measures generated by the spirals belong locally to $H^{-1}(\R^2)$. We say that a measure $\Gamma$ belongs locally to $H^{-1}$ if this measure restricted to any ball $B(0,R_0)$ belongs to $H^{-1}$. The first object to consider
is the Prandtl spiral, whose parametrization in terms of total vorticity is known. Situation seems to be more complicated in the case of Kaden's and Pullin's  spirals with poor analytical description. It appears that also in this case Theorem \ref{Tw.1} is sufficient to prove the desired fact -- it is enough to assume that a spiral satisfies so called Prandtl's similitude laws (see \cite{anton}) to infer $H^{-1}$ local bounds for the measure generated by the spiral from Theorem \ref{Tw.1}. We comment on this issue more in the second part of Section \ref{section3}.

Let us introduce the Prandtl spiral. To describe a position in the plane we use complex variable $z$. By the Prandtl spiral we understand (see \cite{kambe}) the curve given by 
\begin{equation}\label{z_P}
z_P(\Gamma,t) = \left\{\begin{array}{ccl}
\tau^q\Gamma^\nu& \mbox{if} &\Gamma>0,\\ 
-\tau^q|\Gamma|^\nu,& \mbox{if} &\Gamma<0,\end{array}\right.
\end{equation} 
where $\Gamma$ denotes the cumulated strength of the vortex sheet.
 The parameters $\nu=1/2+ib$, $q=1/2+ib\mu$ are complex numbers with $b,\mu$ being real constants. We define $\tau=pt+1$, where $t>0$ denotes time, while $p\in \R$ (in particular we allow $p<0$).

\begin{prop}\label{PrandtlH-1}
For $t\neq -\frac{1}{p}$ the measure generated by the Prandtl spiral belongs locally to $H^{-1}(\R^2)$.
\end{prop}

\proof
The Prandtl spiral has two branches, one is a support of the positive part of the measure defined by the vorticity, another of the negative part. Hence, by Remark \ref{uwaga1}, in order to prove the claim, we consider the measures generated by the positive ($\Gamma>0$) and negative ($\Gamma<0$) branches in \eqref{z_P} separately. Since the proof is the same in both cases, we deal with the positive part only. 

Denote $\Gamma_+(B(0,r))$ the measure of the arc of the positive part of the Prandtl spiral contained in the ball of radius $r$ centered at $0$. We represent $z_P$ in radial coordinates at a fixed time $t$, then 
\[
r=|z_P|=\Gamma^{1/2}\tau^{1/2},
\]
hence
\[
\Gamma_+(B(0,r))=\tau^{-1}r^2 = \frac{1}{tp+1}r^2.
\]
As $\Gamma_+$ is a locally finite Borel measure and thus a Radon measure (\cite[Theorem 2.18]{rudin}) the assumptions of Theorem \ref{Tw.1} are satisfied and we obtain the claim.
\qed
 
From \eqref{z_P} we see that for $t=-\frac{1}{p}$ all the measure of vorticity connected to the Prandtl spiral is contained in one point. 
\begin{rem}\label{uwaga3.1}
For $t= -\frac{1}{p}$ Prandtl's spiral collapses to a point and becomes a measure supported at one point.
\end{rem}
Notice that according to Remark \ref{uwaga3.1} Prandtl spiral is a very interesting example from the point of view of defining global solutions to the 2d Euler equation when initial data is signed measure with nontrivial negative and positive parts. We see that the local kinetic energy, initially being finite, becomes infinite at the time $t=-\frac{1}{p}$. However, in order to use the above property to point the counterexample to global-in-time Delort's solutions one would have to give the meaning to the Prandtl spiral as a solution to the 2d Euler equation. This is still open. 

We would also like to emphasize that when one expects the Prandtl spiral to be a solution of \eqref{B-R}, one meets troubles too, see \cite{kambe}, \cite{saffman}, due to the divergent contributions from the vorticity at infinity. 

Next we would like to analyze further examples of spirals of vorticity known in physics. As those objects do not possess a precise analytical description we will present the methods physicists use to deal with them. Let us notice that it turns out that such objects are expected to satisfy the so-called Prandtl's laws of similitude, which as a consequence yield condition \eqref{warunek} with proper values of $\alpha$, for instance $\alpha=\frac{1}{2}$ in the case of Kaden's spirals. Hence, if the existence of such objects could be proven in a rigorous way, together with the fact that they undergo Prandtl's laws of similitude, the Theorem \ref{Tw.1} would have stated that they do belong locally to $H^{-1}$.

In order to introduce Kaden, Pullin and Anton spirals we need to introduce the method of complex potentials.

First we present the approach based on the complex potential methods. Basic assumptions are as before, the flow is irrotational and incompressible off the interface curve. As in the case of the Prandtl spiral, the complex coordinates are used to describe a point in the plane. If we denote a complex velocity $u:=u_1-iu_2$, then besides an interface curve there exists a complex potential $w$ such that $w'(z)=u(z)$. Actually, $w(x,y)=\phi(x,y)+i\psi(x,y)$, where $\phi$ and $\psi$ are respectively potential and the stream function related to the velocity field $(u_1,u_2)$. Next, one looks for the potential of the velocity of the vortex sheet in the selfsimilar way, namely 
\begin{equation}\label{pot}
w(x,y)=At^{2m-1}\vartheta\left(\frac{x}{t^m},\frac{y}{t^m}\right), 
\end{equation}
where $A$ is a contant, $m>0$ and $\vartheta$ is some profile. The strength of the sheet, $\gamma$, is then given by the jump of the potential across the interface $\gamma:=[w]$.   

An alternative approach consists of introducing the selfsimilar ansatz for both, the position of the curve and the cumulative vortex strength at the point (see \cite{saffman}) in a way
\begin{equation}\label{saf}
z(t, \theta)=t^{m}f(\theta)e^{i\theta}, \;\;\;\Gamma(t,\theta)=t^{2m-1}g(\theta),
\end{equation}
where $f$ and $g$ are selfsimilarity profiles, $m>0$. 

The Kaden spirals (see \cite{kaden}) is a flow related to the choice of potential $w$ with $m=2/3$ above, the profiles $\vartheta$ in \eqref{pot} or $f,g$ in \eqref{saf} being determined for instance by plugging the ansatz in \eqref{B-R}. However there is no rigorous proof that such objects satisfy the 2d Euler equation. Anton (\cite{anton}) argues that the special type of the self-similar flow describing the impulsive flow past a flat plate, also can be found as a potential one of the form \eqref{pot} with $m=2/3$. Once more the argument is informal. More examples of the flows with potentials of the form \eqref{pot}, $0<m<2$, have been studied numerically by Pullin, see \cite{pullin}. Actually, he plugs the form of $z$ and $\Gamma$ as in \eqref{saf} into the Birkhoff-Rott equation  \eqref{B-R}, in turn he arrives at the integral equation which must be satisfied by the profiles $f$ and $g$. Next, he solves the occuring equation numerically. 

The most important information from the point of view of Theorem \ref{Tw.1} is that according to the Prandtl theory, all such flows are supposed to satisfy the similitude laws, see \cite[Section 3]{anton}. Consequently each of them is supposed to satisfy the following relation
\begin{equation}\label{similitude}
\Gamma(B(0,r))=c_2(t)r^{\lambda},
\end{equation}      
where $\lambda>0$ is related to $m$ as $\lambda=2-\frac{1}{m}$. Hence, in the case of the Kaden and Anton spiral $\lambda=1/2$.
Notice that this parameter fits into our Theorem \ref{Tw.1} with $\alpha<1$. We have the following remark.
\begin{rem}\label{last}
Kaden spirals (as well as the Anton and Pullin ones), if they exist and satisfy \eqref{similitude}, by Theorem \ref{Tw.1}, are locally elements of $H^{-1}(\R^2)$. 
\end{rem}

\end{document}